\documentclass{amsart}
\theoremstyle{plain}
\newtheorem{theorem}{Theorem}[section]
\newtheorem{lemma}[theorem]{Lemma}
\newtheorem{proposition}[theorem]{Proposition}
\newtheorem{corollary}[theorem]{Corollary}

\theoremstyle{definition}
\newtheorem{definition}[theorem]{Definition}

\theoremstyle{remark}
\newtheorem{remark}[theorem]{Remark}

\newtheorem{notations}[theorem]{Notations}

\newcommand{\ic}{\ensuremath{\mathcal{I}}}
\newcommand{\gc}{\ensuremath{\mathcal{G}}}
\newcommand{\oc}{\ensuremath{\mathcal{O}}}

\newcommand{\fc}{\ensuremath{\mathcal{F}}}

\newcommand{\jc}{\ensuremath{\mathcal{J}}}

\newcommand{\Pt}{\mathbb{P}^3}

\newcommand{\Pq}{\mathbb{P}^4}
\newcommand{\Pcq}{\mathbb{P}^5}
\newcommand{\Psx}{\mathbb{P}^6}

\newcommand{\Pn}{\mathbb{P}^n}

\newcommand{\bZ}{\mathbb{Z}}

\newcommand{\bQ}{\mathbb{Q}}
\newcommand{\bF}{\mathbb{F}}

\begin{document}

\title[Smooth divisors of projective hypersurfaces.]{Smooth divisors of projective hypersurfaces.}

\author{Ellia Ph.}
\address{Dipartimento di Matematica, via Machiavelli 35, 44100 Ferrara (Italy)}
\email{phe@dns.unife.it}

\author{Franco D.}
\address{Dipartimento di Matematica e Applicazioni "R. Caccioppoli", Univ. Napoli "Federico II", Ple Tecchio 80, 80125 Napoli (Italy)}
\email{davide.franco@unina.it}

\author{Gruson L.}
\address{D\'epartement de Math\'ematiques, Universit\'e de Versailles-St Quentin, Versailles (France)}
\email{gruson@math.uvsq.fr}

\date{12/07/2006}

\maketitle

\section*{Introduction.}

We work over an algebraically closed field of arbitrary characteristic.\\
 Ellingsrud-Peskine (\cite{EP}) proved that smooth surfaces in $\Pq$ are subject to strong limitations. Their whole argument is derived from the fact that the sectional genus of surfaces of degree $d$ lying on a hypersurface of degree $s$ varies in an interval of length $\frac{d(s-1)^2}{2s}$. The aim of the present paper is to show that for smooth codimension two subvarieties of $\Pn$, $n \geq 5$, one can get a similar result with an interval whose length depends only on $s$. The main point is Lemma \ref{TheLemma} whose proof is a direct application of the positivity of $N_X(-1)$ (where $N_X$ is the normal bundle of $X$ in $\Pn$). We get a series of $(n-3)$ inequalities; the first one of which being in \cite{EP}, the second  was obtained in a preliminary version (\cite{EF2}) by an essentially equivalent but more geometric argument.\\
Then we first derive two consequences:\\
1) roughly speaking, (Thm. \ref{Pcq}, Remark \ref{rmkPcq}) the family of "biliaison classes" of smooth subvarieties of $\Pcq$ lying on a hypersurface of degree $s$ is limited\\
2) the family of smooth codimension two subvarieties of $\Psx$ lying on a hypersurface of degree $s$ is limited (Thm. \ref{SpecThm}).\\
The result quoted in 1) is not effective, but 2) is.\\
In the last section we try to obtain precise inequalities connecting the usual numerical invariants of a smooth subcanonical subvariety  $X$ of $\Pn$, $n \geq 5$ (the degree $d$, the integer $e$ such that $\omega _X\simeq \oc _X(e)$, the least degree, $s$, of an hypersurface containing $X$). In particular we prove (Thm. \ref{s=n}): $s \geq n+1$.


\section{Positivity Lemma and some consequences.}

\begin{lemma}[positivity lemma]
\label{TheLemma}
Let $F$ be a rank two vector bundle on a smooth connected variety $X$ of dimension $m$ and let $L$ be an invertible sheaf
such that $h^0(E\otimes L)\not =0$.
Put
$$\frac{1+c_1(L)\,t}{1-c_1 (F)\, t+c_2(F) \, t^2}=\sum u_i t^i $$
in $A_*(X)[[t]]$, where $A_*(X)$ is the Chow ring of $X$ and $t$ is an indeterminate. Assume $F$ is globally generated. Then the $u_i$'s can be represented by pseudo-effective cycles (see \cite{Rob}, 2.2.B), in particular $u_m$ has non-negative degree.
\end{lemma}

\begin{proof}
Set $Q:=P(F)$ (in Grothendieck notation $\mathbf{Proj}(SymF)$) and denote by $p: Q\to X$ the projection. The Chow ring of $Q$ is 
$$
\frac{A_*(X)[x]}{(x^2-c_1(F)\, x+c_2(F))}
$$
(where the indeterminate $x$ corresponds to the tautological quotient 
of $p^*(F)$) and the Gysin map  $p_*: A_*(Q) \to A_*(X)$ sends $\alpha + \beta \, x$ to $\beta$. By hypothesis, there is an effective divisor $D$ of first Chern class  $x+c_1(L)$. Since $F$ is globally generated $x$ is nef and $D.x^i$ is pseudo-effective. Then $p_*(D.x^i)=p_*(x^{i+1}+c_1(L)x^i)=u_i$ (by the formula giving the Gysin map), so $u_i$ is pseudo-effective.  \end{proof}

We will apply the lemma in the following situation:

$X$ is a subvariety of codimension two of $\Pn
 $  (i.e. $n=m+2$) and $F=N_X(-1)$. One knows that $F$ is globally generated because it is a quotient of $T_{\Pn}(-1)$, which is globally generated on $\Pn$. Then we will consider two cases separately:
\begin{enumerate}
\label{2cas}
\item $n=5$.
\item $\omega_X=\oc _X(e)$ for some integer $e$ (by \cite{Barth}  this is always satisfied if $n\geq 6$).
\end{enumerate}

Recall that in the last situation we have an exact sequence
$$
0 \to \oc \to E \to \ic _X (e+n+1) \to 0
$$
where $E$ is a rank $2$ vector bundle on $\Pn $ with Chern classes $c_1(E)=e+n+1$, $c_2(E)=deg(X)$, and that $N_X=E\otimes \oc _X$.

\begin{lemma}
\label{c2Pcq}
Let $X\subset \Pcq $ be a smooth codimension two subvariety of degree d lying on a hypersurface $\Sigma $ of degree
$s$. Denote by $\pi $ the sectional genus of $X$ and assume $X\not \subset Sing \Sigma$. Then one has
$$
0\leq \mu := d(s^2-4s +d)-s(2\pi -2)\leq s(s-1)^3.
$$
\end{lemma}

\begin{proof}
The computations are made in $Num(X)=A(X)/(\text{numerical equivalence})$ (so $Num^3(X)\simeq \bZ$). We denote by $C_i$ the Chern classes of $N_X(-1)$, by $h$ (resp. $k$) the class of $\oc _X(1)$ (resp. $\omega _X$). Finally $y$ will denote the element $c_2(N_X(-s)) \in Num^2(X)$.

By \ref{TheLemma} the $u_i$'s are pseudo-effective. We have $u_2=(s-1)hC_1-C_2=(s-1)^2h^2-y$ and
$$u_3=(s-1)h(C_1^2-C_2)-C_1C_2=(s-1)^3h^3-((s-1)h+C_1)y$$
We know that: $C_1=4h+k$ and $y=(s^2-6s+d)h^2-shk$ (this follows expressing this $c_2$ in function of $c_2(N_X)$ which is $dh^2$ by the self intersection formula). The relation $u_3\geq 0$ is equivalent (in $\bZ$) to:
$$(1)\,\,\,0 \leq (s-1)^3h^3-[(s+3)h+k].[(s^2-6s+d)h^2-shk]=$$
$$-[d(s+3)-21s+1]h^3-(d-9s)h^2k+shk^2$$
Let's write $h^2k$ in function of $\mu$:
$$(2)\,\,\,h^2k=\frac{(s^2-6s+d)h^3-\mu}{s}$$
Apply "Hodge index" to the hyperplane section of $X$:\\ setting
$\delta = (h^2k)^2-h^3.(hk^2)\geq 0$, we get (with $d=h^3$):
$$(3)\,\,\,hk^2=\frac{(h^2k)^2}{d}-\frac{\delta}{d}$$
We partially eliminate $h^2k$:
$$(4)\,\,\,0 \leq -d[d(s+3)-21s+1]+h^2k[-(d-9s)+\frac{s}{d}((s^2-6s+d)\frac{d}{s}-\frac{\mu}{s})]-\frac{s\delta}{d}$$
this yields:
$$(5)\,\,\,0 \leq -d[d(s+3)-21s+d]+h^2k[s(s+3)-\frac{\mu}{d}]-\frac{s\delta}{d}$$
We eliminate $h^2k$:
$$(6)\,\,\,0 \leq -d[d(s+3)-21s+d]+[(s^2-6s+d)d-\mu](s+3-\frac{\mu}{ds})-\frac{s\delta}{d}$$
this can be written
$$(7)\,\,\,0 \leq (s-1)^3d-\mu (2s-3+\frac{d}{s})+\frac{\mu ^2}{ds}-\frac{s\delta}{d}$$
Multiply by $s/d$:
$$(8)\,\,\,0 \leq s(s-1)^3 -\mu [1+\frac{s(2s-3)}{d}] + \frac{\mu ^2}{d^2}- \delta \frac{s^2}{d^2}$$
The relation $u_2 \geq 0$ implies: $\mu \leq d(s-1)^2$, so:
$$(9)\,\,\,0 \leq s(s-1)^3-\mu +\frac{\mu}{d}[-s(2s-3)+\frac{\mu}{d}]-\delta \frac{s^2}{d^2}=$$
$$=s(s-1)^3-\mu +\frac{\mu}{d}[(s-1)^2-s(2s-3)]-\frac{\mu}{d}[(s-1)^2-\frac{\mu}{d}]-\delta \frac{s^2}{d^2}$$
Finally:
$$0 \leq s(s-1)^3 - \mu - \frac{\mu}{d}(s^2-s-1) -\frac{\mu}{d}[(s-1)^2- \frac{\mu}{d}]- \delta \frac{s^2}{d^2}$$
and the lemma follows.

\end{proof}

The last lemma will be used in section \ref{limit}.

In the second case let $s= min\{ t: \, \, h^0(\ic _X(t)\not =0 \}$ and $q=min(s,e+n)$ and notice that
$h^0(N^*_X(q))>0$. Apply the positivity lemma with $L=\oc_(q-e-n)$. Then the $u_i$ can be computed in 
$A_*(\Pn)=\frac{\mathbb{Z}[t]}{t^{n+1}}$ (by abuse of notation we consider $u_i$ as an integer instead of an element of $\mathbb{Z} t^i$) and the positivity lemma applied to $X\cap \mathbb{P}^{i+2}$ says that  $u_i\geq 0$ for $i\leq n-2$.
Let $s_i$ be the Segre classes of $E(-1)$. One has $u_i=c_1(L).s_{i-1}+s_i$. If $s\geq e+n$ one has $L\simeq \oc$, $u_i=s_i$; this case is not new (\cite{HS}), so we focus on the other case ($q=s$).

\begin{remark}
\label{firstui}
The $u_i$ are computed by induction on $i$ by $u_0=1$, $u_1=s-1$, $u_i=(e+n-1)u_{i-1}+(d-e-n)u_{i-2}$.\\
Set $z:=d-s(e+n+1)+s^2$, then the first $u_i$'s are:\\
$u_2=(s-1)^2-z$, $u_3=(s-1)^3-z(e+n+s-2)$ and: $u_4=u_2^2-z(e+n-1)^2$.
\end{remark}

As an immediate consequence  we have:
\medskip

\begin{theorem}[Speciality theorem]
\label{SpecThm}
Let $X \subset \Pn$, $n \geq 5$, be a smooth subvariety of codimension two with $\omega _X \simeq \oc _X(e)$. Let $\Sigma \subset \Pn$ denote an hypersurface of degree $s$ containing $X$. If $X$ is not a complete intersection then:\\
(i) If $n=5$: $e \leq \frac{(s-1)^3}{4}-3-s$ and $d \leq \frac{s(s-1)[(s-1)^2-4]}{4}+1$\\
(ii) If $n \geq 6$: $e \leq \frac{(s-1)^2-n+1}{\sqrt{n-1}}-n+1$ and $d \leq \frac{s[(s-1)^2-n+1]}{\sqrt{n-1}}+1$.
\end{theorem}

\begin{proof}
(i) By $u_3 \geq 0$: $(s-1)^3 \geq z(e+n+s-2)$. Observe that, since $X$ is not a complete intersection, $z=c_2(E(-e-n-1+s)$ is the degree of a codimension two subscheme which is not a complete intersection. By \cite{Ran}, $z \geq n-1$. It follows that: $(s-1)^3 \geq (n-1)(e+n+s-2)$ which gives the bound on $e$. By $u_2 \geq 0$: $d \leq s(n-1+e)+1$ and this gives the bound on $d$.\\
(ii) The proof is similar using $u_4 \geq 0$ instead of $u_3 \geq 0$.
\end{proof}


\section{Application to the  biliaison classes of codimension two subvarieties of $\Pcq$.}
\label{limit}

We recall that a family $\Phi$ of coherent sheaves over an algebraic variety $S$ is \textit{limited} if
there exists an algebraic variety $T$ and a coherent sheaf $\fc $ over $T\times S$ such that for
any member $\gc$ of $\Phi$ there exists a geometric point $t\in T$ such that $\gc $ is isomorphic to the fiber
$\fc _t$ of $\fc $ over $t$.

\begin{theorem}
\label{Pcq}
Fix an integer $s>0$. The family of  sheaves $\ic_{X, \Sigma}([\frac{d}{s}])$, where
\begin{itemize}
\item $\Sigma $ is any integral hypersurface of degree $s$ in $\Pcq$,
\item £d£ is any integer and $X$ is a smooth threefold of degree $d$ lying on $\Sigma$,
\end{itemize}
is limited.
\end{theorem}

\begin{remark}
The corresponding statement for $\Pn$ is
\begin{itemize}
\item false for $n=3$ (for $\Sigma=\mathbb{P}^1\times \mathbb{P}^1$, a quadric in $\Pt$, one gets the sheaves
$\oc(a, -a)$ if $d$ is even),
\item unknown for $n=4$,
\item superseded by the speciality theorem (\ref{SpecThm}) for $n\geq 6$.
\end{itemize}
\end{remark}

Since the degree of $X$ is bounded when $X\subset Sing \Sigma$, the family of the sheaves $\ic_{X, \Sigma}$, with
$X\subset Sing \Sigma$, is clearly limited. Hence in the following we will assume $X\not\subset Sing \Sigma$.

Let $C$ (resp. $S$) denote the intersection of $X$ (resp. $\Sigma$) with a general $\Pt$ in $\Pcq$.

\begin{lemma}
\label{P3}
The sheaves $\ic _{C,S}([\frac{d}{s}])$ form a limited family.
\end{lemma}

\begin{proof}
Due to the existence of Grothendieck $Quot$ scheme, it suffices to show that:
\begin{enumerate}
\item the Hilbert polynomials of these sheaves constitute a finite set
\item there exists an integer $N$ depending only on $s$ such that $\ic _{X,S}([\frac{d}{s}]+N)$ is Castelnuovo-regular.
\end{enumerate}
1) By a direct computation we have:
$$\chi (\ic _{C,S}([\frac{d}{s}])=\binom{s+\epsilon}{3}-\binom{\epsilon}{3}-\frac{\mu}{2s}$$
where $\mu = d(s^2-4s+d)-s(2\pi -2)$, $\epsilon = \frac{d}{s}-[\frac{d}{s}]$. (If $s$ divides $d$, just compare $\chi (\ic _{C,S}(d/s))$ with $\chi (\ic _{\Gamma ,S}(d/s))$ where $\Gamma$ is the complete intersection of $S$ with a surface of degree $d/s$). We conclude with Lemma \ref{c2Pcq}.
\medskip

2) We set $\ic _{C,S}([\frac{d}{s}])=:\fc$ and notice that, for degree reasons, $\fc (s-1)\otimes \oc _{H}$ is Castelnuovo-regular for $H$ a general plane in $\Pt$. Also (since $h^0(\fc (-1))=0$) we have:
$$h^0(\fc (s-1)) \leq \sum _{k=0}^{s-1}h^0(\fc (k) \otimes \oc _H) \leq \sum _{k=0}^{s-1}(sk+1)$$
i.e. $h^0(\fc (s-1))$ is bounded uniformly in $s$. It follows that $h^1(\fc (s-1))$ is bounded uniformly in $s$ (since $h^0$ and $\chi$ are and $h^2(\fc (s-1))=0$), say by $M$. By a classical argument the $h^1$ is strictly decreasing after the regularity of the general plane section (\cite{Sz}) and we deduce: $h^1(\fc (s-1+M))=0$, so $\fc$ is $(s+M)$-regular.
\end{proof}

\begin{lemma}
\label{Phi}
Let $\Phi$ be a family of sheaves on $\Pn$ with the following properties:
\begin{enumerate}
\item  any $\fc \in \Phi $ is locally of depth $\geq2$;
\item for a general hyperplane $H\subset \Pn $ the family of the restrictions of the members of $\Phi $ is limited;
\item $h^0(\fc )$ is bounded uniformly in $\fc \in \Phi$.
\end{enumerate}
Then $\Phi $ is limited.
\end{lemma}

\begin{proof}
By the second assumption we know that the set of the Hilbert polynomials of $\fc \mid _H$ ($\fc \in \Phi $) is finite, so it will be sufficient to prove the following

\textit{Claim}  $h^1(\fc ) $ is bounded uniformly in $\fc \in \Phi$.
\par\noindent
In fact, by assumption $2.$ we know that the 
$h^i(\fc ) $  are bounded uniformly in $\fc \in \Phi$ when $i\geq 2$ because the inequality $h^i(\fc ) \leq \sum_{k\geq 0} h^{i-1}\fc (k)\mid _H$. So, from $1.$ , $3.$ and our claim, it follows that $\mid \chi (\fc) \mid$ is bounded uniformly in $\fc \in \Phi$. So the Hilbert polynomial $P_{\fc} $ of $\fc $ is such that $P_{\fc}(0)$ and ($P_{\fc}(x+1)-P_{\fc}(x)$) form a finite set ($\fc \in \Phi$), this implies that the set $\{P_{\fc}: \, \fc \in \Phi \}$ is finite. A uniform bound on the regularity of $\fc $ is obtained exactly as in the previous lemma.

To prove the claim we look at the exact sequences 
$$
H^0(\fc \mid _H(-k)) \to H^1(\fc (-k-1)) \to H^1(\fc (-k)) \to H^1(\fc \mid _H(-k))
$$
There is an integer $k_0$ independent of $\fc $ such that $h^0(\fc \mid _H(-k_0))=0=h^1(\fc \mid _H(-k_0))$. Since $\fc $ is locally of depth 
$\geq 2$ we also know $H^1(\fc (-k))=0$ for $k>>0$, and so for $k \geq k_0$ by using the above exact sequence.
Then we have $h^1(\fc )\leq \sum _0^{k_0}h^1(\fc \mid _H(-i))$.
\end{proof}

\begin{proof}[Proof of \ref{Pcq}]
By lemma \ref{P3} we know that the family of sheaves $\ic_{X, \Sigma}([\frac{d}{s}]\mid_{\Pt})$ is limited for a general $\Pt \subset \Pcq$. We conclude applying two times lemma \ref{Phi}.
\end{proof}

\begin{remark}
\label{rmkPcq}
\begin{enumerate}
\item If we consider the class of ideals $\ic_{X, \Sigma}$ (as in the theorem) modulo the equivalence
relation identifying two sheaves $\ic $, $\jc $ if $\ic $ is isomorphic to some twist of $\jc $, we could call
them \lq \lq biliaison classes\rq \rq (on a specified hypersurface): if $\ic _{X,\Sigma}\sim \ic_{X', \Sigma'}$
then $\Sigma =\Sigma'$ and $X'$ and $X$ can be linked in $\Sigma $ to the same variety. Then (roughly speaking) the theorem
says that when the degree of the specified hypersurface remains bounded, the set of the corresponding biliaison classes is limited.
\item In contrast with the case $n\geq 6$, we notice that for any $s\geq 2$ one can find ACM, non complete intersection varieties of arbitrary large degree lying on a hypersurface of degree $s$.
\end{enumerate}
\end{remark}

\begin{corollary}
[compare with \cite{BOSS}] The family of smooth threefold in $\Pcq $ which are not of general type is limited.
\end{corollary}
\begin{proof}
According to \cite{BOSS} (proof of Thm. 4.3) we may restrict to the threefolds lying on a hypersurface of degree $12$, so we may fix $s$.
Consider the corresponding family of sheaves $\fc =\ic _{X,\Sigma}([\frac{d}{s}])$, as in Theorem \ref{Pcq}. Then $\omega_X$ is a quotient of $\mathcal{H}om(\ic_{X,\Sigma}, \omega_{\Sigma})=\mathcal{H}om(\fc, \omega_{\Sigma}([\frac{d}{s}]))$. 
Since the family $\Phi$ is limited we can find an integer $k$ (independent of $X$) such that 
$\mathcal{H}om(\fc, \omega_{\Sigma}(k))$ is globally generated. So if $X$ is not of general type one must have $[\frac{d}{s}]< k$
hence $d<s(k+1)$.
\end{proof}


\section{Application to subcanonical codimension two subvarieties of $\Pn$, $n \geq 5$.}

\begin{notations}
\label{gen-notations}
We are now in case 2 of Section \ref{2cas}, so $X$ is the zero-locus of a rank two vector bundle, $E$, of Chern classes $(e+n+1,d)$. For sake of simplicity we consider the Chern polynomial  $e(X)=X^2-C_1X+C_2$ of $E(-1)^*$. Let $\Delta = C_1^2-4C_2$ be its discriminant.
We set $\rho = \sqrt{C_2}$ and write $1-C_1X+C_2X^2=1-2\rho.ch\,t.X+\rho ^2X^2$ with the convention that $t>0$ if $\Delta > 0$ and $t=i\theta$, $0<\theta <\pi$ if $\Delta < 0$ (in this way $ch\,t=cos\theta$, $sh\,t=isin\theta$). Then the roots of $X^2-C_1X+C_2$ are $b= \rho e^t$, $a = \rho e^{-t}$. Finally we set $\sigma = \sqrt z$. \\
If $s_k$ is the $k$-th Segre class of $E(-1)$, i.e.
$$\frac{1}{1-C_1X+C_2X^2}= \sum _{k\geq 0}s_kX^k$$
one deduces from $1-C_1X+C_2X^2=(1-\rho e^tX)(1-\rho e^{-t}X)$, after a partial decomposition, the formula $s_k=\rho ^k\,\frac{sh\,(k+1)t}{sh\,t}$ (to be replaced by $\rho ^k(k+1)$ if $\Delta =0$) and: 
$$u_k=\rho ^k[\frac{s-1}{\rho}\frac{sh\,kt}{sh\,t}-\frac{sh\,(k-1)t}{sh\,t}]$$ (to be replaced by $u_k=\rho ^k[k\frac{s-1}{\rho}-(k-1)]$ if $\Delta =0$). 

\end{notations}

\begin{definition}
\label{defAlpha}  
Let $\Delta > 0$ (resp. $\Delta < 0$), the function $f(x)=\frac{sh\,(x+1)t}{sh\,xt}$ (resp. $g(x)=\frac{sin(x+1)\theta}{sinx\theta}$) is strictly decreasing on $]0,+\infty [$ (resp. $]0,\frac{\pi}{\theta}-1[$). Moreover, if $\Delta > 0$, $lim_{x\to+\infty}f(x)=e^t$. Since $\Delta >0$, $E$ is not stable and $2s < e+n+1$, since $0<z=e(s-1)$, we have $s-1 <a$ hence $s-1\leq a=\rho e^{-t}$, $e^t \leq \frac{\rho}{s-1}$. We conclude that there exists an unique $\alpha$ such that $f(\alpha )=\rho /(s-1)$ (resp. such that $g(\alpha )=\rho /(s-1)$). Notice that in the case $\Delta < 0$, $(\alpha +1)\theta < \pi$.\\
Similarly considering the function $v(x)=(x-1)/x$ we define, in case $\Delta =0$, $\alpha$ such that: $\frac{\alpha}{\alpha +1}=\frac{s-1}{\rho}$. By Lemma \ref{TheLemma}: $\alpha \geq n-3$.
\end{definition}

\begin{lemma}
\label{panoplie}
With notations as above we have:\\
if $\Delta > 0$: 
$$\frac{\sigma}{sh\, t} = \frac{s-1}{sh\, \alpha t} = \frac{\rho}{sh\, (\alpha +1)t} = \frac{e+n-s}{sh\, (\alpha +2)t} $$
if $\Delta < 0$:
$$\frac{\sigma}{sin\, \theta} = \frac{s-1}{sin\, \alpha \theta } = \frac{\rho}{sin\, (\alpha +1)\theta } = \frac{e+n-s}{sin\, (\alpha +2)\theta } $$
if $\Delta =0$:
$$\alpha + 2 = \frac{e+n-s}{\sigma} = \frac{\rho}{\sigma}+1 = \frac{s-1}{\sigma}+2 $$
\end{lemma}

\begin{proof}
First assume $\Delta > 0$. By definition $z=e(s-1)$. Inserting $s-1=\frac{\rho sh\,\alpha t}{sh\,(\alpha +1)t}$, we get 
$$z=\rho ^2[\frac{sh^2(\alpha t)}{sh^2(\alpha +1)t}-2ch\,t.\frac{sh(\alpha t)}{sh(\alpha +1)t} +1]=\rho ^2[\frac{sh^2(\alpha t)}{sh^2(\alpha +1)t} - \frac{sh(\alpha -1)t}{sh(\alpha +1)t}]$$. For the last equality check that $sh\,(\alpha +1)t+sh\,(\alpha -1)t=2ch\, t.sh\,(\alpha t)$. Finally: $\rho ^2[\frac{sh^2(\alpha t)}{sh^2(\alpha +1)t} - \frac{sh(\alpha -1)t}{sh(\alpha +1)t}] = [\frac{\rho sh\,t}{sh\,(\alpha +1)t}]^2$. For this check that $sh^2(\alpha t)-sh(\alpha-1)t.sh(\alpha +1)t=sh^2(t)$. We conclude that $\frac{\sigma}{sh\,t}=\frac{\rho}{sh(\alpha +1)t}$. This proves the first three equalities. For the last one: 
$$\frac{\rho}{sh\,(\alpha +1)t}=\frac{s-1}{sh\,\alpha t}=\frac{2\rho ch\,t\,-(s-1)}{2sh\,(\alpha +1)t\,ch\,t - sh\,\alpha t}$$
To conclude observe that: $2\rho ch\,t = e+n-1$ and $2sh\,(\alpha +1)t\,ch\,t - sh\,\alpha t = sh\,(\alpha +2)t$.\\  
The proof in case $\Delta < 0$ is similar. If $\Delta =0$, observe that $z = e(s)=(s-\rho -1)^2$, hence 
$\sigma = \rho -s+1$.
\end{proof}

\begin{remark}
Observe that when $\Delta < 0$ and $s=e+n$, then $sin(\alpha +2)\theta =0$.
\end{remark}

\begin{proposition}
Keeping notations as above, we have, for $n \geq 5$:
$$ e+n-s \leq (n-1)^{-\frac{1}{n-4}}(s-1)^{\frac{n-2}{n-4}}
$$
and
$$
d\leq s [1+(n-1)^{-\frac{1}{n-4}}(s-1)^{\frac{n-2}{n-4}} ]
$$
\end{proposition}

\begin{proof}
First of all we assume $\Delta >0$
 and we observe that $f(t)=logsh\,t$ is concave. Since 
$\alpha t=\frac{1}{\alpha}t + \frac{\alpha -1}{\alpha}(\alpha+1)t$ 
we have $f(at)=f(\frac{1}{\alpha}t + \frac{\alpha -1}{\alpha}(\alpha+1)t)\geq 
\frac{1}{\alpha}f(t)+\frac{\alpha -1}{\alpha}f((\alpha+1)t)$. Taking the exponentials we find
$$
sh\, \alpha t\geq (sh\,t)^{\frac{1}{\alpha}}(sh(\alpha +1)t)^{\frac{\alpha -1}{\alpha}} \,  \, \, (+)
$$
Similarly, writing $\alpha t=\frac{2}{\alpha +1}t + \frac{\alpha -1}{\alpha +1}(\alpha+2)t$ and exponentiating the
inequality coming from the concavity of $f(t)$ we get
$$
sh\, \alpha t\geq (sh\,t)^{\frac{2}{\alpha +1}}(sh(\alpha +2)t)^{\frac{\alpha -1}{\alpha +1}} \, \, \, (++)
$$
By \ref{panoplie} 
$\frac{\sigma}{sh\, t} = \frac{h-1}{sh\, \alpha t} = \frac{\rho}{sh\, (\alpha +1)t}= \frac{e+n-h}{sh\, (\alpha +2)t} $,
hence $(+)$ gives
$$
s-1\geq (\sigma )^{\frac{1}{\alpha}}(\rho)^{\frac{\alpha -1}{\alpha}} 
$$
and $(++)$ gives
$$
s-1 \geq (\sigma )^{\frac{2}{\alpha +1}}(e+n-s)^{\frac{\alpha -1}{\alpha +1}} 
$$
from which follow
$$
\rho \leq \sigma ^{\frac{-1}{\alpha -1}}(s-1)^{\frac{\alpha}{\alpha -1}}, \, \, \, \, \, \, \, \, \, \, \,
e+n-s \leq \sigma ^{\frac{-2}{\alpha -1}}(s-1)^{\frac{\alpha +1}{\alpha -1}}.
$$
and finally $d=\rho ^2 +e+n \leq s( 1+ \sigma ^{\frac{-2}{\alpha -1}}(s-1)^{\frac{\alpha +1}{\alpha -1}}).$

In order to conclude the case $\Delta >0$ it suffices to show that 

$\sigma ^{\frac{-2}{\alpha -1 }}(s-1)^{\frac{\alpha +1}{\alpha -1}}\leq (n-1)^{\frac{-1}{n-4}}(s-1)^{\frac{n-2}{n-4}}$. Since $z=\sigma ^2 \geq n-1 $ we have 
$\frac{(s-1)^2}{\sigma ^2}\leq \frac{(s-1)^2}{n-1}$ so
$\sigma ^{\frac{-2}{\alpha -1 }}(s-1)^{\frac{\alpha +1}{\alpha -1}}=(\frac{(s-1)^2}{\sigma ^2})^{\frac{1}{\alpha -1}}(s-1)
\leq (\frac{(s-1)^2}{n-1})^{\frac{1}{\alpha -1}}(s-1)$ and we are done because $\alpha \geq n-3$.

The case $\Delta < 0$ ($\Delta =0$) can be proved the same way by using $f(t)=logsin\, t$ ($f(t)=log\,t$) which is concave as well for $t\in ]0, \pi [$.
\end{proof}

In some sense the next proposition improves Theorem \ref{SpecThm}, except in the case $\Delta >0$ where the bound depends on $\Delta$, hence on $e$.

\begin{proposition}
\label{d<Ms2}
Let $X \subset \Pn$, $n \geq 4$ be a smooth codimension two subvariety with $\omega _X\simeq \oc _X(e)$. If $X$ is not a complete intersection, then: 
\begin{enumerate}
\item If $\Delta >0$, then $d < M^2s^2+sM\sqrt \Delta$, where $M=\frac{n-2}{n-3}$.
\item If $\Delta \leq 0$, then $d < M^2s(s-1)+s$
\end{enumerate} 
\end{proposition}

\begin{proof}
1) By \ref{panoplie} 
$$\frac{\rho}{s-1}=\frac{sh(\alpha+1)t}{sh\alpha t}$$
and 
$$\frac{sh(\alpha+1)t}{sh\alpha t}\leq \frac{sh(n-2)t}{sh(n-3) t}$$ since $\alpha \geq n-3$ (\ref{defAlpha}).
One can check that 
$$\frac{sh(n-2)t}{sh(n-3) t}\leq e^t \frac{n-2}{n-3}$$ so 
$$\frac{\rho}{s-1}\leq e^t \frac{n-2}{n-3}$$ and 
$$a=\frac{\rho}{e^t}\leq (s-1)\frac{n-2}{n-3}$$
Then we have 
$$w:= a- (s-1)\leq (s-1)[\frac{n-2}{n-3}-1]=\frac{s-1}{n-3}< \frac{s}{n-3}$$ and
$$d=(a+1)(b+1)=(s+w)^2+ \sqrt\Delta (s+w).$$ Finally we get
$d< (\frac{n-2}{n-3})^2s^2 + s \sqrt \Delta \frac{n-2}{n-3}$.
\\
2) Assume first $\Delta < 0$. By Lemma \ref{panoplie} we have: $\frac{\rho}{s-1}=\frac{sin((\alpha +1) \theta )}{sin((\alpha )\theta)} \leq \frac{\alpha + 1}{\alpha}$, indeed $\frac{sin x}{x}$ is decreasing on $0 < x \leq \pi$. It follows that $\rho \leq M(s-1)$. Since $\rho = \sqrt{ab}=\sqrt{d-e-n}$, we get the result taking into account the inequality: $s(e+n+1-s)\leq d$ ($z \geq 0$). The case $\Delta =0$ follows directly from $u_k\geq 0$ (see \ref{gen-notations}), taking into account the inequality $s(e+n+1-s) \leq d$.  
\end{proof}

\begin{remark}
\label{rmkd<Ms2}
Observe the limiting ($n \to +\infty$) case of Prop. \ref{d<Ms2}, 1): $d \leq s^2 +s\sqrt{\Delta}$, which can occur only for $X$ a complete intersection $(a+1,b+1)$.
\end{remark}
The aim of the remaining of the paper is to improve the bound $s \geq n-1$ of \cite{Ran} (resp. $s \geq n$ if $5\leq n\leq 6$, \cite{EF}). We will distinguish several cases according to the sign of the discriminant, $\Delta$, of $X$.

\begin{proposition}
\label{sDelta>0}
Let $X \subset \Pn$, $n \geq 4$, be a smooth subvariety of codimension two. Assume $X$ is not a complete intersection.
\begin{enumerate}
\item if $\Delta \geq 0$, then $s-1 \geq (n-3)\sqrt{n-1}$ and $e \geq (2n-4)\sqrt{n-1}-n$.
\item if $\Delta < 0$ and $e+n+1 \geq 2s$, then $s-1 \geq \frac{2}{\pi}(n-3)\sqrt{n-1}$ and\,\,\, $e \geq \frac{2}{\pi}(2n-4)\sqrt{n-1}-n$.
\end{enumerate}
\end{proposition}

\begin{proof}
1.) Assume first $\Delta > 0$. By Lemma \ref{panoplie}: $\frac{s-1}{\sigma} = \frac{sh(\alpha t)}{sh\,t}\geq \alpha$. Since $\sigma \geq \sqrt{n-1}$ and $\alpha \geq n-3$, we get the result. In the same way, from Lemma \ref{panoplie}: 
$\frac{e+n-s}{\sigma}=\frac{sh(\alpha +2)}{sh\,t}\geq \alpha +2 \geq n-1$, hence: $e+n-s \geq (n-1)\sqrt{n-1}$ and the result follows.

Assume now $\Delta = 0$.
We have $z=(s-a-1)^2$. Since $z \geq n-1$, it follows that $a+1-s \geq \sqrt {n-1}$ (note that $a+1 > s$ if $X$ is not a c.i.), so $a \geq s-1+\sqrt {n-1}$ and we get
$s-1 \geq \frac{(n-3)}{(n-2)}(s-1+\sqrt {n-1})$ hence $s-1 \geq (n-3)\sqrt{n-1}$. We conclude as above since
$\frac{e+n-s}{\sigma}= \alpha +2 \geq n-1$.\\
2.) By Lemma \ref{panoplie} $\frac{e+n-s}{s-1 }=\frac{sin(\alpha +2)\theta}{sin\alpha \theta}$. The assumption $e+n+1 \geq 2s$ implies that $sin(\alpha +2)\theta \geq sin\,\alpha \theta$. This in turn, implies $\alpha \theta < \frac{\pi}{2}$ (we have $(\alpha +2)\theta < \frac{3\pi}{2}$, cf Definition \ref{defAlpha}). By Lemma \ref{panoplie} $\frac{s-1}{\sigma}=\frac{sin(\alpha \theta)}{sin\, \theta}$. Since $sin\alpha x /sin\,x$ is decreasing, we get: $\frac{s-1}{\sigma}>\frac{1}{sin(\pi/2\alpha)}$, hence $s-1 > \frac{2}{\pi}(n-3)\sqrt{n-1}$. The proof for $e$ is similar using $\frac{e+n-s}{\sigma}=\frac{sin(\alpha +2)\theta}{sin\, \theta}$ of Lemma \ref{panoplie}.

\end{proof}

\begin{lemma}
\label{s-DeltaNeg-L2}
Let $X \subset \Pn$, $n \geq 6$, be a smooth subvariety of codimension two. Assume $\Delta <0$ and $e+n+1-2s \leq 0$.\\
(i) If $n \geq 6$, then $s \geq n+2$\\
(ii) If $n \geq 8$, then $s \geq 3n/2$
\end{lemma}

\begin{proof}
(i) If $n \geq 6$, then $e \geq n+2$ (\cite{HS}), hence $(n+2)+n+1 \leq e+n+1\leq 2s$, thus $s \geq n+2$.\\
(ii) As above using $e \geq 2n-1$ (\cite{HS} Cor.3.4 (i)).
\end{proof}

\begin{remark}
\label{mn}
In case $\Delta < 0$, we may proceed as follows: by Lemma \ref{panoplie}: $s-1 = \rho \frac{sin(\alpha \theta)}{sin(\alpha +1)\theta}$, so $s-1= \delta \frac{sin(\alpha \theta)}{sin\,\theta .sin(\alpha +1)\theta}$, where $\delta = \frac{1}{2}\sqrt{-\Delta}=\rho sin\,\theta$, $\theta (\alpha +1)<\pi$. Let's denote by $m(\alpha )$ the minimum of the function $\varphi (\theta )=\frac{sin(\alpha \theta)}{sin\,\theta.sin(\alpha +1)\theta}$ on $]0,\frac{\pi}{\alpha +1}[$. This minimum is reached for the solution, $\beta$, of $\frac{sin(\alpha +1)\beta}{sin\,\beta}=\sqrt{\alpha +1}$ and is an increasing function of $\alpha$. So we have: $s-1 \geq \delta m(\alpha) \geq \delta m(n-3) \geq \frac{\sqrt{-\Delta _{min}(n)}}{2}.m(n-3)$, where $-\Delta _{min}(n)$ is the minimal value of $-\Delta$ allowed by the Schwarzenberger conditions on $\Pn$ (see \cite{HS}). It is possible to compute an approximated value of $m_n:=m(n-3)$. For instance we have: $m_5=1,6949$, $m_6=2,2845$, $m_7=2,8203$, $m_8=3,3233$ (and $m_{40}=16,1647$). Since $-\Delta _{min}(8)=119$, we get: $s-1 \geq 19$ if $n=8$, which is better than $12=\frac{3.8}{2}$.\\
Let $E$ be a rank two vector bundle on $\Pn$ with Chern classes $c_1,c_2$ (and $\Delta = c_1^2-4c_2$ not a square). Let $R = \bZ[X]/(X^2-c_1X+c_2)$. The Schwarzenberger condition says that: $Tr_{\bQ R/\bQ}\binom{\xi +k}{n} \in \bZ$ for $\xi$= class of $X$, $k \in \bZ$. Let $p$ be a prime number, then we have three cases:
\begin{enumerate}
\item inert ($pR$ is prime)
\item decomposable ($R/pR \simeq \bF _p\times \bF _p$)
\item ramified ($p\,|\,\Delta$)
\end{enumerate}
\textit{Claim: If there exists a rank $2$ vector bundle $E$ of Chern classes $(c_1,c_2)$ on $\Pn$, then for each prime 
$p<n$, the discriminant $\Delta =c_1^2-4c_2$ is a square mod $p$ (possibly $0$).}
\begin{proof} Assuming the contrary, we may suppose that $n-1=p$ is a prime such that $\Delta $ is not a square mod $p$.
 Let $\xi $ be a root of $X^2-c_1X+c_2$ in $\bF _{p^2}$. Then $Tr_{\bF _{p^2}/\bF _p}(\xi(\xi+1)\dots (\xi +p))\equiv -\Delta $ $mod\, p$, since $(\xi+1)\dots (\xi +p)=F(\xi)-\xi$  where $F$ is the Frobenius automorphism of 
$\bF _{p^2}$ and $Tr_{\bF _{p^2}/\bF _p}(\xi (F(\xi)-\xi))\equiv -\Delta$. So, if $x$ is the image of $X$ in $R$, one has 
$Tr_{R/\bZ}(x(x+1)\dots (x+p))$ is not divisible by $p$ and $Tr_{\bQ R/\bQ } \binom{x+p}{p+1}=\chi (E)$ is not an integer.
Contradiction.
\end{proof}
By \cite{Serre} pp. 134-135 one knows that there exists some prime $p<n$ such that $\Delta$ is not a square \textit{mod} $p$,
when $n\geq c(log|\Delta |)^2$ under Generalized Riemann Hypothesis or when
$n\geq 2 (|\Delta |)^A$ without restrictions. This means that $|\Delta _{min} (n)|\geq e^{\sqrt\frac{n}{c}}$ 
under GRH or $|\Delta _{min} (n)|\geq (\frac{n}{2})^{\frac{1}{A}}$. 
\end{remark}


\begin{corollary}
\label{3n/2}
Let $X \subset \Pn$, $n \geq 11$, be a smooth codimension two subvariety. If $X$ is not a complete intersection, then $s \geq \frac{3n}{2}$.
\end{corollary}

\begin{proof}
Assume first $\Delta < 0$. If $e+n+1 > 2s$, then, by \ref{sDelta>0} it is enough to check that $1+\frac{2}{\pi}(n-3)\sqrt{n-1} \geq \frac{3n}{2}$ if $n \geq 11$. If $e+n+1 \leq 2s$, then $s \geq \frac{3n}{2}$ by Lemma \ref{s-DeltaNeg-L2}.\\
If $\Delta \geq 0$, by Proposition \ref{sDelta>0}, it is enough to check that $1+(n-3)\sqrt{n-1}\geq \frac{3n}{2}$ if $n \geq 11$ (actually this is satisfied for $n \geq 8$).
\end{proof}

\begin{theorem}
\label{s=n}
Let $X \subset \Pn$, $n \geq 4$, be a smooth codimension two subvariety. Assume $ch(k)=0$. If $n < 6$ assume $X$ subcanonical. If $h^0(\ic _X(n))\neq 0$, then $X$ is a complete intersection.
\end{theorem}

\begin{proof}
For the case $n=4$ we refer to \cite{EFG}.  If $n=5$ by \cite{EF} we may assume $s=5$ and by \cite{BC} $e \geq 3$. From $u_3\geq 0$ (see Remark \ref{firstui}) we get: $z \leq \frac{(s-1)^3}{e+n+s-2}$, i.e. $z \leq 5$. In fact $4\leq z \leq 5$, since $z \geq n-1$ (\cite{Ran}). Arguing as in \cite{EF} Lemma 2.6, every irreducible component of $Z_{red}$ appears with multiplicity, hence $Z$ is either a multiplicity $z$ structure on a linear subspace or is contained in a cubic hypersurface. The last case is not possible (\cite{Ran}). In the first case by \cite{Ma} (or also \cite{Ma2} observing that the proof of the main theorem works in the case of a codimension two linear subspace of $\Pcq$), $Z$ is a complete intersection.\\
If $6 \leq n \leq 7$. By $u_4 \geq 0$: $f(z)=\frac{[(s-1)^2-z]^2}{z} \geq (e+n-1)^2$. Since $f(z)$ is decreasing and $z \geq n-1$, $f(n-1)\geq (e+n-1)^2$, i.e. $(s-1)^2 \geq \sqrt {n-1}(e+n-1)+n-1$. By \cite{HS}: $e \geq n+2$, so $(s-1)^2 \geq \sqrt {n-1}(2n+1)+n-1$, but this inequality is not satisfied if $s \leq n$, $6 \leq n \leq 7$.\\
Now assume $8 \leq n \leq 10$. If $\Delta \geq 0$, we conclude with Proposition \ref{sDelta>0}. If $\Delta < 0$ we conclude  by Remark \ref{mn} ($s \geq 20$ if $n=8$).
For $n \geq 11$, we conclude with Corollary \ref{3n/2}.
\end{proof}


\end{document}